\documentclass[12pt]{article}
\newtheorem{lemma}{Lemma}
\newtheorem{thm}{Theorem}
\usepackage{latexsym}
\usepackage{amssymb}
\usepackage{amsmath}
\begin{document}
\thispagestyle{empty}

\vskip 20pt
\begin{center}
{\Large New Bounds on van der Waerden-type Numbers for Generalized 3-term Arithmetic Progressions}
\vskip 15pt
{\bf Patrick Allen}\\
{\it Department of Mathematics}\\
{\it University of Maryland - Baltimore County}\\
\vskip 10pt
{\bf Bruce M. Landman}\\
{\it Department of Mathematics}\\
{\it University of West Georgia}\\
\vskip 10pt
{\bf Holly Meeks}\\
{\it Department of Mathematics}\\
{\it University of West Georgia}
\end{center}
\vskip 30pt
\begin{abstract}
Let $a$ and $b$ be positive integers with $a \leq b$. An $(a,b)$-triple
is a set $\{x, ax+d, bx+2d\}$, where $x,d \geq 1$. Define
$T(a,b;r)$ to be the least positive integer $n$ such that any
$r$-coloring of $\{1,2,\ldots,n\}$ contains a monochromatic $(a,b)$-triple. Earlier results gave an upper bound on $T(a,b;2)$ that is a fourth degree polynomial in $b$ and $a$, and a quadratic lower bound. A new upper bound for
$T(a,b;2)$ is given that is a quadratic. Additionally, lower bounds are given for the case in which $a=b$, updated tables are provided, and open questions are presented.
\end{abstract}

\vspace{.3in}


\vspace{.04in}


\section{Introduction}
Van der Waerden [10] proved that for any positive integers $n$ and $r$, there is a least positive integer $w(n;r)$ such that every $r$-coloring of $\{1,2,\ldots,w(n;r)\}$ admits a monochromatic $n$-term arithmetic progression.
The estimation of the function $w(n;r)$ is a notoriously difficult problem of Ramsey theory.
Numerous results analogous to van der Waerden's theorem, where the family of arithmetic progressions is replaced by other families, have been considered. Ramsey properties for a variety of such families,  along with a summary of what is known about $w(n;r)$, may be found in [6]. Other recent results of this type may be found in [1], [4], [7], and [8].

In [5], a statement analogous to van der Waerden's theorem was considered for a certain generalization of 3-term arithmetic progressions. Specifically,
for given integers $a$ and $b$ with $a \leq b$, any ordered triple of the form ${x, ax+d,
bx+2d}$, where $x$ and $d$ are positive integers, is called an {\em $(a,b)$-triple}.  We see that the (1,1)-triples are just the 3-term arithmetic progressions. Similar to the function $w(k;r)$, for each pair $(a,b)$, $1 \leq a \leq b$, and $r$ a positive integer, define $T(a,b;r)$ to be the least positive integer $t$ (if it exists) such that
every $r$-coloring of $\{1,2,\ldots,t\}$ contains a monochromatic $(a,b)$-triple.
Hence, $T(1,1;r)$ has the same meaning as $w(3;r)$.

Along with van der Waerden's theorem, another classical result of Ramsey theory is due to Rado [9], which deals with the existence of monochromatic solutions of systems of linear equations under finite colorings. Much is known about the Rado numbers corresponding to a single linear homogeneous equation. That is, given $c,d,e \in \mathbf{Z}$, the {\em Rado  number}
$R(cx + dy = ez;r)$ is the least  positive  integer $n$ such that every $r$-coloring of $\{1,2,\ldots,n\}$ contains a monochromatic integer solution to $cx + dy = ez$ (note that $x$ and $y$ need not be distinct).
We observe that $(x,y,z)$ is an $(a,b)$-triple if and only if
  $z = 2y+(b-2a)x$ and $y > ax$, where $x$ is any positive integer. Thus, $T(a,b;r)$ may be considered to be a ``restricted" Rado number.

As most of the work in this paper deals with the situation in which $r=2$, for convenience we will denote $T(a,b;r)$ more simply as $T(a,b)$. In [5] it was shown that $T(a,b)$ exists if and only if $b \neq 2a$. Further, the following upper and lower bounds were found for $T(a,b)$:

\begin{thm} Let $1 \leq a \leq b$.
\begin{enumerate}
\item If $b > 2a$, then
\[ 2b^{2}+5b - 2a + 4 \leq T(a,b) \leq 4a(b^{3}+b^{2}-3b-3)+2b^{3}+4b^{2}+6b. \]
\item If $b < 2a$, then
\[ 3b^{2} - 2ab + 5b - 2a +4 \leq T(a,b) \leq 4a(b^{3}+2b^{2}+2b)-4b^{2} .\]
\end{enumerate}
\end{thm}

For  the special case in which $b=a$, a quadratic upper bound on $T(a,b)$ was given; and for $b=2a-1$, a cubic upper bound was provided.

The main result of this paper, given in Section 2, shows that there is a quadratic upper bound on $T(a,b)$, for all $(a,b)$ such that $b \neq 2a$.

For a positive $r$, we say that $(a,b)$ is {\em $r$-regular} if $T(a,b;r)$ exists. By van der Waerden's theorem, $(1,1)$ is $r$-regular for all $r$. In [2], [3], and [5], results on the {\em degree of regularity} of $(a,b)$-triples, i.e., the largest $r$ such that $T(a,b)$ exists, was investigated. In particular, in both [2] and [3] it was shown, independently, that (1,1) is the only pair that is $r$-regular for all $r$. In Section 3, we give updated tables on values of $T(a,b)$ and on the degree of regularity of $(a,b)$.

We employ the following additional notation and terminology. We denote by $[1,m]$ the set $\{1,2,\ldots,m\}$. An {\em $(a,b)$-valid} coloring of a set is a coloring that avoids monochromatic $(a,b)$-triples.

\section{An Upper Bound on $T(a,b)$}

In this section we give an improvement over Theorem 1, namely a quadratic upper bound on $T(a,b)$ for all $a,b$ such that $b \neq 2a$.
Our upper bound makes use of the following lemma.

\begin{lemma}
Let $a,b,k,M,i \in {\mathbf Z}^{+}$ and $a \leq b$.
Let $\chi$ be an $(a,b)$-valid 2-coloring of $[1,M]$ with
$\chi(k) \neq \chi(k+i)$.

\begin{enumerate}

\item[(a)] Assume
\begin{enumerate}
\item[(i)] $x > ak + ib/2$
\item[(ii)] $M \geq \max\{x,2(x-ak)+bk\}$
\item[(iii)] $i(b-2a)$ is even
\item[(iv)] $\chi(x) = \chi(k)$.
\end{enumerate}
Then $\chi(x - i(b-2a)/2) = \chi(k).$

\item[(b)] Assume
\begin{enumerate}
\item[(i)] $y > b(k+i)$
\item[(ii)] $M \geq \max\{y,y-i(b-2a)\}$
\item[(iii)] $y$ and $b(k+i)$ have the same parity
\item[(iv)] $\chi(y) = \chi(k+i).$
\end{enumerate}
Then $\chi(y-i(b-2a)) = \chi(k+i)$.
 \end{enumerate}
\end{lemma}
{\em Proof.}
(a) By (i) and (iii), $(k+i,a(k+i)+(y-b(k+i))/2),y)$ is an
$(a,b)$-triple. Hence, since $\chi$ is valid on $[1,M]$, by (iv) we have

\begin{equation}
\chi\left(a(k+i) + \frac{1}{2}(y-b(k+i) \right) = \chi(k) .
\end{equation}

Now, $(k,a(k+i) + (y-b(k+i))/2, y - i(b-2a))$ is an $(a,b)$-triple in $[1,M]$ (with $d = ai + (y-b(k+i))/2$, and thus
(1) implies $\chi(y-i(b-2a)) = \chi(k+i)$.

(b) Since $x > ak$, we know $(k,x,bk+2(x-ak))$ is an $(a,b)$-triple, and hence, by (iv),
\begin{equation}
\chi(bk+2(x-ak)) = \chi(k+i) .
\end{equation}

By (i) and (iii), $(k+i,x-i(b-2a)/2,bk+2(x-ak))$ is an $(a,b)$-triple (with $d = x - ak -ib/2$). By (2) and hypothesis (iv), the result follows.
\hfill
$\Box$

We may now obtain a quadratic upper bound on $T(a,b)$, thereby improving the fourth degree polynomial upper bound of Theorem 1.
\begin{thm}
Let $1 \leq a < b$ with $b \neq 2a$.  Then
\[ T(a,b) \leq \begin{cases}

  7b^2 - 6ab +13b-10a & \text{for $b$ even, $b > 2a$} \\

  14b^{2} - 12ab + 26b - 20a & \text{for $b$ odd, $b > 2a$} \\

  3b^{2} +2ab + 16a & \text{for $b$ even, $b < 2a$} \\

  6b^2+4ab+8b+16a & \text{for $b$ odd, $b < 2a$}. \\

\end{cases}\]
\end{thm}
{\em Proof.}
Case 1: $b-2a$ even, $b > 2a$. Let $M = 7b^2 - 6ab + 13b - 10a$, and let $\chi$ be an arbitrary 2-coloring
of $[1,M]$. We will show that there is a monochromatic $(a,b)$-triple under $\chi$. Since $(1,a+1,b+2)$ is an $(a,b)$-triple, we may assume there exists $c \leq b+1$ with $\chi(c) \neq
\chi(c+1)$.  Let

\[ \beta=2 \left\lfloor \frac{2a(c+1)}{b-2a}\right\rfloor + 2c +2,\]
and let $x_0 = (\beta + c + 1)(b-2a)$.

We consider two subcases:

\noindent
 Subcase (A): $\chi(x_0) = \chi(c)$.  We will apply Lemma 1(a) with $x = x_0$, $k=c$ and $i=1$. Provided the hypotheses of the lemma hold, this will give
 that $\chi(x_0) = \chi(x_{0} -(b-2a)/2)$. We will then repeatedly apply the lemma, in turn, to each of $x_j = x_{0}-j(b-2a)/2$, $j=1,2,\ldots, \beta + 2c +2$, until
 we obtain that $\chi(\frac{\beta}{2}(b-2a)) = \chi(c)$.

 In order to do so, we now check that the hypotheses of Lemma 1(a) hold at each step in this process.
Clearly, hypotheses (iii) and (iv) hold.
For hypothesis (i), we have that, for each $x_j$,

\begin{eqnarray*}
x_j & \geq & \left(\frac{\beta +1}{2}\right)(b-2a)  \\
& = & \left( \left\lfloor \frac{2a(c+1)}{b-2a}\right\rfloor + c + \frac{3}{2} \right) \\
& > & \left(\frac{2a(c+1)}{b-2a} + c + \frac{1}{2}\right)(b-2a)  \\
& = & bc + \frac{b}{2} + a  \\
& > & ac + \frac{b}{2}.
\end{eqnarray*}

Furthermore, hypothesis (ii) holds for each $j$, $0 \leq j \leq \beta + 2c +2$, since
\begin{eqnarray} 2x_j + c(b-2a) & \leq & (2\beta + 3c + 2)(b-2a) \nonumber \\
& = & \left( 4 \left\lfloor \frac{2a(c+1)}{b-2a}\right\rfloor  + 7c + 6 \right)(b-2a) \nonumber \\
& \leq & 8a(c+1) + (7c+6)(b-2a) \nonumber \\
& \leq & 8a(b+2) + (7b+13)(b-2a) \nonumber \\
& = & M.
\end{eqnarray}

As a result, we may conclude that the $(a,b)$-triple $(c,\frac{\beta}{2}(b-2a),(\beta+c)(b-2a))$ is monochromatic.

\vspace{.05in}

\noindent
Subcase (B): $\chi(x_0) = \chi(c+1)$. Taking $k=c$ and $i=1$, we will repeatedly apply Lemma 1(b), beginning with $y = x_0$, until we obtain
$\chi(\frac{\beta}{2}(b-2a) = \chi(c+1)$. This will give the monochromatic $(a,b)$-triple

\[ (c+1, \frac{\beta}{2}(b-2a),x_0) .\]

We now check that, indeed, the hypotheses of Lemma 1(b) are satisfied for each $y = x_0 - j(b-2a)$, $j=0,1, \ldots, \frac{\beta}{2} + c$.
 By assumption, hypothesis (iv) holds, and since $b$ and $\beta$ are both even, we see that hypothesis (iii) holds.
 Now,
 \[ \left( \frac{\beta}{2} + 1 \right) (b-2a) > bc + \frac{b}{2} + a + \frac{b}{2}  -a = b(c+1), \]
 so that hypothesis (i) holds for each application of the lemma. Since $y \leq 2x_{0} + c(b-2a)$,  by (3) we know that hypothesis
 (ii) holds.

In both subcases, $[1,M]$ contains a monochromatic $(a,b)$-triple, completing the proof in this case.

\noindent
Case 2: $b$ odd, $b-2a>0$.

Let $M =14b^{2}-12ab+26b - 20a$. Assume $\chi$ is a 2-coloring of $[1,M]$. We must show that $\chi$ is not valid on $[1,M]$.
Since $(2,2a+2,2b+4)$ is an $(a,b)$-triple, we may assume there is an even $c$, with $c \leq 2b+2$, such that $\chi(c) \neq \chi(c+2)$.
Let
\[ x_0 = (2\beta + c+2)(b-2a),  \mbox{   where   }  \beta = 2 \left\lfloor \frac{ac+2a}{b-2a} \right\rfloor .\]

We consider two subcases.

\noindent
Subcase (A): $\chi(x_0) = \chi(c)$. We shall repeatedly apply Lemma 1(a), with $k = c$ and $i = 2$, beginning with $x = x_0$, to obtain

\begin{equation}
 \chi(c) = \chi(x_0) = \chi(x_0 - (b-2a)) = \cdots = \chi(\frac{\beta}{2}(b-2a)).
\end{equation}

In order to do so, let us check that the hypotheses of the lemma hold at each step in the process. Clearly, (iii) and (iv) are true.

The least value of $x$ for which we employ the lemma is $u = (\frac{\beta}{2} + 1)(b-2a)$. Since

\begin{eqnarray*} u & = & \left(\left\lfloor \frac{ac+2a}{b-2a}\right\rfloor + \frac{c}{2} + 2 \right)(b-2a)  \\
& > & ac + 2a + \left(\frac{c}{2} + 1\right)(b-2a)  \\
& = & b + \frac{b}{2}c > b + ac,
\end{eqnarray*}
hypothesis (i) of Lemma 1(a) holds for each iteration of the lemma.

Finally, to verify that hypothesis (ii) holds at each step, it is sufficient to show that $M \geq 2x_0 + c(b-2a)$. This does hold since

\begin{eqnarray}
2x_0 + c(b-2a) & = & \left(8 \left\lfloor \frac{ac+2a}{b-2a} \right\rfloor + 7c + 12\right)(b-2a) \nonumber \\
& \leq & 8ac + 16a + (7c+12)(b-2a) \nonumber \\
& \leq & 8a(2b+2) + 16a + (14b+26)(b-2a) \nonumber \\
& = & M.
\end{eqnarray}

Hence, from (4), the $(a,b)$-triple $(c,\frac{\beta}{2}(b-2a),(\beta+c)(b-2a))$ is monochromatic.

\noindent
Subcase (B): $\chi(x_0) = \chi(c+2)$. Repeatedly applying Lemma 1(b), starting with $y = x_0$, until we obtain $\chi(c+2) = \chi(\beta(b-2a)$,  will yield the monochromatic $(a,b)$-triple $(c+2,\beta(b-2a),(2\beta+c+2)(b-2a))$. The lemma does apply in each instance since,  for each value of $y$, $y$ is even, $M \geq y$ by (5), and
\begin{eqnarray*}
y & \geq & (\beta + 2) (b-2a)  \\
& > & \left( \frac{2a(c+2)}{b-2a} + c + 2\right)(b-2a)  \\
& = & b(c+2).
\end{eqnarray*}

In both subcases, $[1,M]$ has a monochromatic $(a,b)$-triple.

\noindent
Case 3. $b$ even, $b < 2a$.

Let $M = 3b^{2} + 2ab + 16a$, and let $z = \beta(2a-b)$, where
\[ \beta = 2\left\lfloor \frac{ac+b/2}{2a-b} \right\rfloor - c+ 2 .\]

Note that
\begin{equation}
z > \left(\frac{2ac+b}{2a-b} - c \right)(2a-b) = b(c+1).
\end{equation}
As in Case 1, we may assume there is $c \leq b+1$ with $\chi(c) \neq \chi(c+1)$.
We have two subcases.

\noindent
Subcase (A): $\chi(z) = \chi(c)$. By repeatedly applying Lemma 1(a), taking $k=c$ and $i=1$, beginning with $x=z$, we will obtain
\[
 \chi(z) = \chi\left(z+ 1/2(2a-b)\right) = \chi(z+(2a-b)) = \cdots = \chi(z+(\beta+c)(2a-b)).
\]
Hence, the $(a,b)$-triple $(c,z,(2\beta+c)(2a-b))$ is monochromatic. We see from (6) that hypothesis (i) of the lemma
holds in each instance. Obviously, hypotheses (iii) and (iv) also hold. To show that hypothesis (ii) holds, note first that, taking $x$ as in the lemma, the largest value of $x$
to which we need to apply the lemma is $x = z + (\beta + c -1/2)(2a-b).$
Now, for this value of $x$, since $c \leq b+1$, we have
\begin{eqnarray}
2x + c(b-2a) & = & 2(2\beta + c -1/2)(2a-b) + c(b-2a) \nonumber \\
& = & (4\beta + c - 1)(2a-b) \nonumber \\
& \leq & 8ac + 4b + (-3c+7)(2a-b) \nonumber  \\
& \leq & M .
\end{eqnarray}

\noindent
Subcase (B): $\chi(z) = \chi(c+1)$. We shall apply Lemma 1(b), repeatedly, beginning with $y = z$, until we obtain
\[ \chi((2\beta + c +1)(2a-b)) = \chi(c+1). \]
Now, hypothesis (i) of Lemma 1(b) holds by (6), and hypothesis (ii) is immediate from (7), showing that the lemma may be applied at each step. This gives the monochromatic $(a,b)$-triple $(c+1,z,(2\beta + c + 1)(2a-b))$.

In both subcases, there is a monochromatic $(a,b)$-triple in $[1,M]$.

\noindent
Case 4. $b$ odd, $b < 2a$.

Let $M = 6b^{2} + 4ab + 8b+16a$.
Let $z = \beta(2a-b)$, with
\[ \beta = 2\left\lfloor\frac{ac+b}{2a-b}\right\rfloor + 2 -c. \]
Note that

\begin{equation}
\beta(2a-b) > \left(\frac{2ac+2b}{2a-b} - c\right)(2a-b) = b(c+2).
\end{equation}
We may assume there is a $c \leq 2b+2$ such that $\chi(c) \neq \chi(c+2)$, and consider two subcases.

\noindent
Subcase (A): $\chi(z) = \chi(c)$. We repeatedly apply Lemma 1(a), taking $i=2$ and $k=c$, beginning with $x =z$, which yields the monochromatic
$(a,b)$-triple $(c,z,(2\beta+c)(2a-b))$. The lemma is applicable at each step since hypothesis (i) holds by (8), and hypothesis (ii) holds because the largest value $x$ to which
the lemma will be applied is $u=(2\beta+c-1)(2a-b)$ and
\begin{eqnarray}
 2u - c(2a-b)& = & (4\beta+c-2)(2a-b) \nonumber \\
 & \leq & 8ac + 8b + (6-3c)(2a-b) \nonumber \\
 & \leq & (2b+2)(2a+3b) +2b + 12a \nonumber \\
 & = & M.
\end{eqnarray}
\noindent
Subcase(B): $\chi(z) = \chi(c+2)$. Letting $k=c$ and $i=2$, we apply Lemma 1(b), beginning with $y=z$ until we obtain the monochromatic
$(a,b)$-triple $(c+2,z,(2\beta+c+2)(2a-b))$. By assumption, hypothesis (iii) of Lemma 1(b) holds. Hypothesis (i) holds by (8), and hypothesis
(ii) holds by (9).

In either case, there is a monochromatic $(a,b)$-triple in $[1,M]$.
\hfill $\Box$

\noindent
{\bf Remark 1.} Lemma 1 can be extended to include negative $i$, as long as $k > -i$. To do so, we need only to change hypothesis (i) of (a) to
$x > \max\{ak, ak + ib/2\},$ and hypothesis (i) of (b) to $ y > \max\{b(k+i),b(k+i)-2ai\}$. However, applications of this extended lemma do not
yield better bounds than those of Theorem 2.

\section{Using $r$ Colors}

As noted in the introduction, $(1,1)$ is the only pair that is $r$-regular for all $r$. We also noted that the pairs $(a,b)$ for which $b = 2a$ are the only pairs whose degree of regularity is 1. In [2] it is shown that
d.o.r.$(a,b) \leq 5$ for all pairs other than (1,1). It is known that certain infinite families of pairs $(a,b)$ have d.o.r.$(a,b) \leq 4$, others have d.o.r.$(a,b) \leq 3$, and others have d.o.r.$(a,b) = 2$ (see [2],[3]). On the other hand, the only $(a,b)$ for which it has been verified that d.o.r.$(a,b) > 2$ is $(a,b) = (2,2)$; in fact, $T(a,b;3) = 88$ (see [3]).

In the following table, which updates a table in [5], we summarize what is known about the degree of regularity for small values of $a$ and $b$. The values in the table are based on various results and proofs from [2],[3],and [5].

\[
\begin{tabular}{l|cccc|}
$b\backslash a$   &   1   &   2   &   3   &   4 \\ \hline
1&$\infty$&-&-&-\\
2&1&3-4&-&-\\
3&2&2&2-3& - \\
4&2&1&2-3&2-3 \\
5&2-4&2&2&2-3\\
6&2-3&2&1&2\\
7&2-4&2-4&2&2\\
8&2-5&2-3&2&1\\
9&2-5&2-4&2-3&2\\
10&2-3&2-4&2-3&2\\ \hline
\end{tabular}
\]
\begin{center}
{\bf Table 1 Degree of Regularity for Small Values of $a$ and $b$}
\end{center}

In [5], it was shown that $T(a,a;2) \geq a^{2}+3a+4$. The next result improves this slightly, and also provides lower bounds on $T(a,a;3)$ and $T(a,a;4)$.

\begin{thm}
\begin{enumerate}
\item $T(a,a;2) \geq a^{2} + 3a + 8$ for $a \geq 4$.
\item  $T(a,a;3) \geq 3a^{3} + 4a^{2} + 5a + 8$ for $a > 1$.
\item $T(a,a;4) \geq 7a^{4}+12a^{3}+6a^{2}+9a+16$ for $a > 1$.
\end{enumerate}
\end{thm}
{\em Proof.}
 (1) The 2-coloring that is red on $[a^{2}+2,a^{2}+2a+1] \cup \{a^{2}+3a+4, a^{2}3a+6\}$ and blue on $[1,a+1] \cup [a^{2}+2a+2,a^{2}+3a+3] \cup \{a^{2}+3a+5,a^{2}+3a+7\}$ can easily be shown to be valid.

\noindent
(2) Let $M = 3a^{3} + 4a^{2} + 5a +7$. It suffices to show that there exists an $(a,a)$-valid 3-coloring of
$[1,M]$. Consider the coloring defined as follows. Color the intervals $R_{1}= [1,a+1]$ and $R_{2} = [a^{3}+2a^{2}+2a+2,2a^{3}+3a^{2}+3a+3]$ red, the intervals $B_{1} = [a+2,a^{2}+2a+1]$ and $B_{2} = [2a^{3}+3a^{2}+3a+4,3a^{3}+4a^{2}+5a+7]$ blue, and the interval
$G=[a^{2}+2a+2,a^{3}+2a^{2}+2a+1]$ green.

Assume there is a red $(a,a)$-triple. The only $(a,a)$-triple having two terms in $R_{1}$ is $(1,a+1,a+2)$, so at most one term of the red triple is in $R_1$. If the first term, $x$, does belong to $R_1$ and the second term belongs to $R_{2}$, then $d \geq a^{3}+a^{2}+a+2$. But then $ax+2d$ lies outside of $R_2$, a contradiction. Hence, $x \in R_2$; but then, since $a \geq 2$, $ax+2d \geq a^{4}+2a^{3}+2a^{2}+2a+2 > 2a^{3}+3a^{2}+3a+3$, which is not possible.

If there is a blue $(a,a)$-triple, it clear that at most one term can be in $B_{1}$. If the first term, $x$, of the triple is in $B_{1}$, then
$d \geq a^{3}+a^{2}+2a+4$,  and hence $ax+2d$ is not blue. Therefore, $x \in B_{2}$, implying that $ax+2d \not\in B_{2}$, so there is no blue $(a,a)$-triple.

Finally, it is clear that there is no green $(a,a)$-triple, completing the proof of (2).

\noindent
(3) The proof follows the same reasoning as (2) and is left to the reader.
\hfill $\Box$

In [5], a table of values of $T(a,b)$ is given. The following table updates that list. New values are indicated by the symbol *.

\[
\begin{tabular}{l|ccccccc|}
$a\backslash b$ & 1 &2 &3 &4 &5 &6 &7\\ \hline
1&9&$\infty$&39&58&81&$108^{*}$&$139^{*}$ \\
2&&16&46&$\infty$&139&$114^{*}$&$159^{*}$\\
3&&&39&60&114&$\infty$&$247^{*}$\\
4&&&&40&87&$144^{*}$&$\geq214$ \\
5&&&&&70&100&$\geq150$\\
6&&&&&&78&$120^{*}$\\
7&&&&&&&95\\ \hline
\end{tabular}
\]
\begin{center}
{\bf Values of $T(a,b)$}
\end{center}

We conclude with some questions that we find intriguing:
\begin{itemize}
\item Does there exist a pair $(a,b)$ whose degree of regularity is greater than 3?
\item Characterize the pairs $(a,b)$ whose degree of regularity is greater than 2.
\item It is known [5] that \[ T(1,b) \geq 2b^{2} + 5b + 6\] for all $b \geq 3$.
Is this lower bound the actual value of $T(1,b)$ for all $b \geq 3$? From Table 1, we see that this is true for all $b$, $3 \leq b \leq 7$.
\item It was shown in [5] that \[ T(a,2a-1) \geq 16a^{2} - 12a + 6, \] for all $a \geq 2$.
We wonder if this inequality is an equality; the answer is yes for $a=2,3.$
\item Is it true that $T(b,b) \leq T(a,b)$ for all $1 \leq a \leq b$?
\end{itemize}

\end{document}